\documentclass[sn-mathphys-num]{sn-jnl}% Math and Physical Sciences Numbered Reference Style
\usepackage[english]{babel}
\usepackage[utf8]{inputenc}
\usepackage{graphicx}%
\usepackage{multirow}%
\usepackage{amsmath,amssymb,amsfonts}%
\usepackage{amsthm}%
\usepackage{mathrsfs}%
\usepackage[title]{appendix}%
\usepackage{xcolor}%
\usepackage{textcomp}%
\usepackage{manyfoot}%
\usepackage{booktabs}%
\usepackage{algorithm}%
\usepackage{algorithmicx}%
\usepackage{algpseudocode}%
\usepackage{listings}%

\usepackage{bm}
\usepackage{tikz,pgfplots}
\usepackage{tikz-3dplot}
\usetikzlibrary{patterns}
\usetikzlibrary{arrows.meta, 3d}

\theoremstyle{thmstyleone}%
\theoremstyle{thmstyletwo}%

\theoremstyle{thmstylethree}%

\begin{document}

\title[Computational algorithm for downward continuation of gravity anomalies]{Computational algorithm for downward continuation of gravity anomalies\footnote[2]{Preprint submitted to arXiv}}
\author[1]{\fnm{D. K.} \sur{Ivanov}}\email{ivanov.gx@s-vfu.ru}
\author[2]{\fnm{L. N.} \sur{Temirbekova}}\email{laura-nurlan@mail.ru}
\author*[3,4]{\fnm{P. N.} \sur{Vabishchevich}}\email{vab@cs.msu.ru}

\affil[1]{\orgdiv{Yakutsk Branch of the Regional Scientificand Educational Mathematical Center ``Far Eastern Center of Mathematical Research''}, \orgname{Ammosov North-Eastern Federal University},
 \orgaddress{\street{Kulakovskogo 48}, \postcode{677013}, \city{Yakutsk}, \country{Russia}}}
 
\affil[2]{\orgdiv{Faculty of Mathematics, Physics and Informatics}, \orgname{Abai Kazakh National Pedagogical University}, \orgaddress{\street{Dostyk avenue 13}, \postcode{050010}, \city{Almaty}, \country{Kazakhstan}}}

\affil*[3]{\orgdiv{Faculty of Computational Mathematics and Cybernetics}, \orgname{Lomonosov Moscow State University}, \orgaddress{\street{Leninskie Gory 1}, \postcode{119991}, \city{Moscow}, 
\country{Russia}}}

\affil[4]{\orgdiv{Institute of Mathematics and Information Science}, \orgname{Ammosov North-Eastern Federal University},
 \orgaddress{\street{Kulakovskogo 48}, \postcode{677013}, \city{Yakutsk}, \country{Russia}}}
 
\abstract{
The downward continuation of potential fields from the Earth's surface into the subsurface is a critical task in gravity exploration, as it helps to identify the sources of gravity anomalies.
This problem is often addressed by solving a first-kind integral equation using regularization techniques to stabilize an inherently unstable process.
A similar approach is used in our work, where the continued field is represented as the potential of a simple layer or its vertical derivative.
The constancy of the density sign of this equivalent simple layer preserves the sign of anomalies, provided that the layer's surface encloses all anomalous sources.
This constraint is a key feature of our algorithm for the downward continuation of potential fields.
To enforce, for instance, non-negativity in the simple layer density, we employ the NNLS (Non-Negative Least Squares) method.
The efficiency of the proposed method is demonstrated on model examples.
}



\maketitle

\section{Introduction}

An essential class of geophysical problems involves gravimetry, which examines heterogeneities in the Earth's crust and their influence on the gravitational field \cite{gupta2011encyclopedia,mudrecova1990gravi}.
The primary measured quantity is gravitational acceleration, and the gravitational anomaly is defined as deviations from the normal gravity field after applying necessary corrections.
Analysis of gravity survey data provides insights into mass distribution in the upper crust
\cite{Parker1994,blokh2009interp}.

The direct gravimetric problem is associated with the calculation of the gravitational field generated by a mass with known physical and geometric parameters.
In inverse gravimetric problems these parameters are deduced from observed field data \cite{Crossley2013}, typically collected the Earth's surface, e.g. on earth's surface, airborne surveys and satellite missions.
Solutions of direct and inverse problems are used in data processing, qualitative interpretation of the observed data.
Among these tasks, transformation, upward and downward continuation of the gravitational field, and characterization of the shape and density of anomalies are highlighted separately.
When interpreting observed gravity fields, solving the forward problem makes it possible to predict and compare gravity anomalies for geological structures.
Numerical methods for solving direct and inverse problems of geophysics specifically gravimetry and magnitometry are well-studied
\cite{Vogel1992,zhdanov2002geophysical,freeden2018handbook,yagola2014obratnye,eppelbaum2019geophysical}.

The gravitational potential is defined as the Newtonian (volume) potential.
The measured gravitational field is associated with the derivatives of the first- and higher-order potentials \cite{blakely1996potential}.
The observations themselves are carried out on the Earth's surface, taking into account the elevation \cite{dolgal2022gravimetry} and the Earth's sphericity for large-scale surveys.
There are two main continuation tasks with respect to the observation surface and source locations: upward continuation and downward continuation \cite{pilkington2017potential}.
The gravitational potential and its derivatives outside the anomalies satisfy the Laplace equation.
Thus, we have problems of continuing solutions of elliptic equations beyond the boundary of the computational domain, which reduces to Cauchy problems for elliptic equations.
Field downward continuation (towards sources) allows for clearer localization of the anomalies \cite{pavsteka2012regcont}.

The Cauchy problem for elliptic equations belongs to the class of ill-posed problems \cite{lavrentiev2006operator}, where small perturbations of the input data (measurements at the boundary) amplify solution errors.
Approximate solutions employ regularization techniques \cite{tikhonov1977solutions,kabanikhin2011inverse}.
For conditionally well-posed evolutionary problems, various approaches of approximate solution have been developed.
For instance, we highlight \cite{samarskii2007numerical} the Tikhonov regularization method, methods of equation perturbation (the quasi-inversion method), and separately note the algorithms with perturbation of initial conditions (non-local boundary value problems).
In modern computing practice, iterative methods are most widely used \cite{alifanov1995extrem,hansen1998rank,vasin2013ill,Bakushinsky2011}, where the number of iterations acts as a regularization parameter.

For geopotential field continuation from one layer (irregular surface or level plane) to another, potentials are often approximated by equivalent sources, while the original generating sources are unknown.
For brevity, we mention works \cite{pilkington2017potential,dolgal2022history} and references therein.
For these problems, field-based and source-based methods can be distinguished.
Field-based approaches are primarily based on simple and/or double-layer potentials derived from Green's identities.
For approximate solutions of Cauchy problems for elliptic equations in the case of unbounded domains, integral equations are more effective with well-developed numerical methods \cite{Hackbusch1995,kythe2011computational}.
Such an approach was applied to potential downward continuation problems by \cite{tikhonov1968prodolgenie}.
In the approximate solution of ill-posed problems, a fundamental increase in accuracy is achieved \cite{tikhonov1995numerical} by taking into account additional information about the solution.
In a series of works \cite{strakhov2002s,stepanova2019approx,stepanova2020separation}, the method of S-approximations and its modifications for separating geopotential fields by the sum of potentials of simple and double layers was proposed.

When continuing potential fields from one plane to another, the problem can be more effectively solved by Fourier transforms in the spectral domain \cite{pavsteka2012regcont}.
Wavelet transformation techniques are used in multilevel source approximation in \cite{li2010rapid,balk2016effective}; in \cite{balk2016effective}, this approach with Haar basis functions is compared with the quadtree algorithm.
Wavelet analysis using the Legendre basis is applied in \cite{temirbekov2022numerical} for the solution of the first-kind integral equation.

Source-based methods \cite{dampney1969equivalent} operate with elementary homogeneous objects such as points, prisms, etc., which generate parameterized potentials.
Among recent works, in \cite{leonov2024solving}, the authors pointed out the instability of source approximation, which can be effectively solved by a new approach based on improving the condition number of a matrix with the help of the minimal pseudoinverse matrix method.
To continue the field into the lower half-space, the problem becomes unstable, requiring stable techniques based on iterative processes or regularization \cite{zeng2013adaptive,li2023stable}.

In \cite{vabishchevich2024computational}, a two-dimensional problem of potential field continuation towards anomaly sources was considered using the integral equation method.
The two-dimensional problem corresponds to masses extended along the horizontal direction, and gravity is calculated using the logarithmic potential.
In this work, we extend the approach to the case of three-dimensional problems.
The gravitational field is represented as a simple-layer potential or its derivative.
It is shown that the equivalent layer density is positive (negative) for positive (negative) density anomalies if the surface encloses all anomalies.
The continuation process is solved using the well-known method of non-negative least squares (NNLS) \cite{lawson1995solving}, enforcing the required sign of the simple-layer densities.

The article is structured as follows.
The article begins by introducing the relevance of gravimetry and the challenges associated with downward continuation of gravity anomalies (Section 1). Section 2 formulates the mathematical framework, presenting integral equations for potential field continuation. Section 3 derives constraints on the equivalent layer density to enforce non-negativity, while Section 4 details the NNLS-based methodology for solving the ill-posed inverse problem. Section 5 validates the approach through 3D simulations, demonstrating robustness against noise and accurate depth localization. Finally, Section 6 synthesizes key findings, highlighting the algorithm’s effectiveness in estimating anomaly sources.  

\section{Problem statement}
\begin{figure}[htb]
	\centering
	\includegraphics[width=0.6\textwidth]{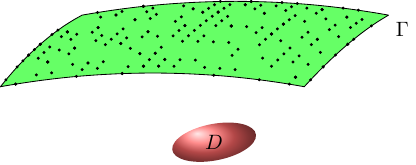}
	\caption{Anomalous geodensity body $D$ and observation surface $\Gamma$ with irregular grid marked by dots.}
	\label{f:1}
\end{figure}

We consider the problem of continuing gravitational fields measured on the Earth's surface.
Gravitational anomalies correspond to density contrasts between subsurface masses and the surrounding rock formations.
Let us define the anomalous density distribution as $\rho(\bm x)$, $\bm x = (x_1, x_2, x_3)$.
We assume that the anomaly occupies a domain $D \subset \mathbb{R}^3$, and outside this domain the density is considered to be zero.
Observation data are generally provided on an irregular surface $\Gamma$, as illustrated in Figure~\ref{f:1}.
We assume that $\Gamma \subset \mathbb{R}^3$ and $\Gamma \cap D = \varnothing$, representing the Earth's surface with relief taken into account.

The gravitational field due to the anomalous mass is determined by the volume potential:
\begin{equation}\label{e:1}
	u(\bm x) = \int_D K(\bm x, \bm y) \rho(\bm y) \, d\bm y, \quad \bm x \in \Gamma,
\end{equation}
where the kernel is given by
\[
	K(\bm x, \bm y) = G \frac{1}{|\bm x - \bm y|}, \quad |\bm x|^2 = \sum_{i=1}^{3} x_i^2,
\]
and $G$ is the universal gravitational constant.

The vertical derivative of the gravitational potential is a commonly observed quantity during surface surveys:
\begin{equation}\label{e:2}
\frac{\partial u}{\partial l} = g(\bm x), \quad \bm x \in \Gamma,
\end{equation}
where $\partial/\partial l$ denotes the vertical derivative. For local and medium-scale surveys, the direction $l$ coincides with the $x_3$-axis, which is oriented upwards in our case.
In inverse gravimetry problems, the goal is to obtain information about the spatial domain $D$ containing anomalies, the density distribution $\rho(\bm x)$ for $\bm x \in D$, and the gravitational field near the anomalies, based on measurements of $g(\bm x)$ on $\Gamma$.

\begin{figure}[htb]
	\centering
	\includegraphics[width=0.8\textwidth]{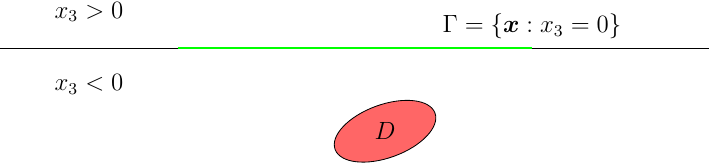}
	\caption{Anomalous geodensity body $D$ and a horizontal observation plane $\Gamma$.}
	\label{f:2}
\end{figure}

For simplicity, we consider a horizontal observation plane $\Gamma = \{\bm x \in \mathbb{R}^3: x_3 = 0\}$, as shown in Figure~\ref{f:2}.
The continuation of the gravitational field into the upper half-space ($x_3 > 0$) reduces to a boundary value problem for the Laplace equation:
\[
\Delta u = 0, \quad x_3 \ge 0
\]
with the boundary condition on $\Gamma$:
\[
\frac{\partial u}{\partial l} (\bm x) = g(\bm x), \quad x_3 = 0.
\]
Additionally, the potential $u$ satisfies the decay condition at infinity:
\[
\lim_{|\bm x| \rightarrow \infty} u(\bm x) = 0, \quad \bm x \in \mathbb{R}^3.
\]

The continuation of the gravitational field toward the source region, in the absence of information about its geometry and density, is a classical ill-posed problem.
To obtain a meaningful solution, additional constraints on the solution space must be imposed.
The continuation problem in the lower half-space is formulated as follows: find $u(\bm x)$ for $x_3 < 0$ such that
\[
\Delta u = 0, \quad \bm x \notin D.
\]

Gravitational field continuation can be formulated as a Cauchy boundary value problem for an elliptic equation, for which various computational approaches have been studied \cite{samarskii2007numerical}.
Among these, we highlight the Tikhonov regularization method based on a variational formulation interpreted as an optimal control problem.
The quasi-inversion method \cite{lattes1969method} involves solving a perturbed boundary value problem, namely with perturbations in the equation and boundary/initial conditions (non-local boundary value problems).
Iterative methods are widely used in the numerical solution of ill-posed and inverse problems.
In \cite{vab2002num}, an iterative method is proposed based on successive improvements of the continued field by solving standard boundary value problems.

\begin{figure}[t!]
	\centering
	\includegraphics[width=0.6\textwidth]{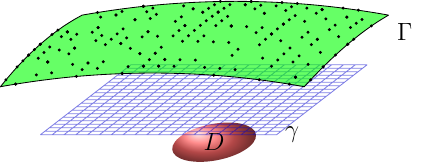}
	\caption{Anomalous geodensity body $D$, observation surface $\Gamma$, and continuation surface $\gamma$.}
	\label{f:3}
\end{figure}

We approximate the gravitational potential~\eqref{e:1} and its corresponding field~\eqref{e:2} using a simple-layer potential.
The approximation is constructed on an arbitrary surface $\gamma$ located just below the observation surface $\Gamma$.
The potential is then approximated as
\[
u(\bm x) = \int_\gamma \mu(\bm y) K(\bm x, \bm y) \, ds(\bm y), \quad \bm x \in \Gamma,
\]
where $\mu(\bm y)$ is the surface density function.
The corresponding integral equation for $\mu(\bm x)$ is formulated as
\[
\int_\gamma \mu(\bm y) \tilde K(\bm x, \bm y) \, d\bm y = g(\bm x), \quad \bm x \in \Gamma,
\]
with known observations $g(\bm x)$ on $\Gamma$.
Here, the kernel $\tilde K$ is obtained by differentiating $K$ along the direction $l$:
\[
\tilde{K}(\bm x, \bm y) = \frac{\partial K}{\partial l}(\bm x, \bm y).
\]
The continuation surface $\gamma$ may be irregular but sufficiently smooth; in Figure~\ref{f:3}, it is shown as a horizontal plane.

\section{Class of positive solutions}
\begin{figure}[htb]
	\centering
	\includegraphics[width=0.6\textwidth]{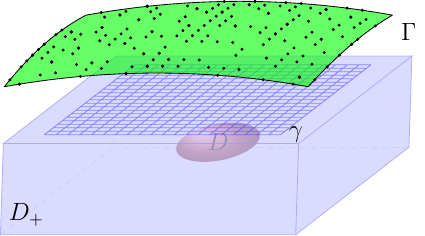}
	\caption{Extended domain $D^+$ enclosing the anomalous geodensity body $D$. Observation surface $\Gamma$ and continuation plane $\gamma$ are shown as parts of the boundary of $D^+$.}
	\label{f:4}
\end{figure}

Consider an auxiliary boundary value problem in an extended domain $D^+$ such that it contains the anomalous body $D \subset D^+$ and $D^+ \cap \Gamma = \varnothing$, see Figure~\ref{f:4}.
Let the function $v(\bm x)$, $\bm x \in D^+$, satisfy the Poisson equation, analogous to the volume potential $u(\bm x)$:
\begin{equation}\label{e:3}
 \Delta v = - \chi_D(\bm x) \rho(\bm x), \quad \bm x \in D^+,
\end{equation}
where the characteristic function of the domain $D$ is defined as
\[
\chi_D(\bm x) = 
\begin{cases}
1, & \bm x \in D,\\
0, & \bm x \notin D.
\end{cases}
\]
A homogeneous Dirichlet condition is imposed on the boundary $\partial D^+$:
\[
v(\bm x) = 0, \quad \bm x \in \partial D^+.
\]

Taking into account both the equation and boundary condition for $v$, we apply the third Green’s identity to $v$ with the fundamental solution $K$ of the Laplace equation. This yields the following integral representation outside the domain $D^+$:
\[
 \int_{D^+} \chi_D(\bm y) \rho(\bm y) K(\bm x, \bm y) \, d\bm y = - \int_{\partial D^+} \frac{\partial v}{\partial n} (\bm y) K(\bm x, \bm y) \, ds(\bm y), \quad \bm x \in \Gamma,
\]
where $n$ denotes the outward unit normal to the boundary $\partial D^+$.

The left-hand side of this equation corresponds to the volume potential:
\[
u(\bm x) = - \int_{\partial D^+} \frac{\partial v}{\partial n}(\bm y) \, K(\bm x, \bm y) \, ds(\bm y), \quad \bm x \in \Gamma.
\]
This dimensionality reduction technique is commonly used in forward modeling of potential fields \cite{vabishchevich1983economical}.

Now consider problems where the density distribution $\rho(\bm x)$ has a constant sign within the domain $D^+$.
Without loss of generality, we assume $\rho(\bm x) \ge 0$.
The maximum principle for elliptic equations ensures that $v(\bm x) \ge 0$ for all $\bm x \in D^+$.
This, in turn, implies a constraint on the normal derivative: ${\partial v}/{\partial n} \le 0$ on $\partial D^+$.
We define the surface density $\mu(\bm x)$ as
\[
\mu(\bm x) = - \frac{\partial v}{\partial n}(\bm x), \quad \bm x \in \partial D^+,
\]
which leads to the non-negativity condition:
\[
\mu(\bm x) \ge 0, \quad \bm x \in \partial D^+.
\]
Hence, the gravitational potential can be represented as a simple-layer potential:
\[
u(\bm x) = \int_{\partial D^+} \mu(\bm y) K(\bm x, \bm y) \, ds(\bm y), \quad \bm x \in \Gamma.
\]

If we consider $D^+$ as the lower half-space, then the boundary $\partial D^+$ becomes an infinitely extended horizontal plane that separates the observation surface $\Gamma$ from the gravitating masses $D$.
From a computational standpoint, we restrict this infinite boundary to a finite continuation surface $\gamma \subset \partial D^+$, located directly beneath $\Gamma$; see Figure~\ref{f:4}.
Neglecting the influence of distant sources, the gravitational potential is approximated by the simple-layer potential on $\gamma$:
\[
u(\bm x) = \int_\gamma \mu(\bm y) K(\bm x, \bm y) \, ds(\bm y).
\]
The inverse problem is then reformulated as follows: find the surface density function $\mu(\bm x)$ satisfying the integral equation
\begin{equation}\label{e:4}
\int_\gamma \mu(\bm y) \tilde{K}(\bm x, \bm y) \, ds(\bm y) = g(\bm x), \quad \bm x \in \Gamma,
\end{equation}
subject to the non-negativity constraint
\begin{equation} \label{e:5}
\mu(\bm x) \ge 0, \qquad \bm x \in \gamma.
\end{equation}
This restriction to the class of non-negative solutions is crucial for developing a stable computational algorithm for solving the inverse problem~\eqref{e:4}--\eqref{e:5}.

\section{Computational algorithm}

Let $w_\Gamma$ be a finite set of observation points on the surface $\Gamma$:
\[
w_\Gamma = \{\bm x^{(i)} \mid \bm x^{(i)} = (x_1^{(i)}, x_2^{(i)}, 0), \ i = 0, 1, \ldots, N\},
\]
where the measurements $g(\bm x^{(i)})$ are taken.
For the horizontal continuation plane $\gamma$ at depth $h$, we construct a rectangle grid:
\[
w_\gamma = \{\bm y^{(j)} \mid \bm y^{(j)} = (y_1^j, y_2^j, h), \ j = 0, 1, \ldots, M\},
\]
with nodes $\bm y^{(j)}$ representing the centers of surface cells. The depth value $h$ is user-defined and can be freely chosen.

Let the measured quantity be the vertical derivative of the gravitational potential, i.e., $\partial/\partial l = -\partial/\partial x_3$ (assuming the $x_3$ axis points upward). In this case, the kernel $\tilde K$ takes the form
\[
\tilde K(\bm x, \bm y) = G \frac{x_3 - y_3}{|\bm x - \bm y|^3},
\]
so that positive mass anomalies produce positive gravity anomalies.

The integral equation~\eqref{e:4} is discretized using the rectangle quadrature rule over the grid $w_\gamma$:
\[
\sum_{j=0}^{M} 
\tilde K(\bm x^{(i)}, \bm y^{(j)}) \, \mu(\bm y^{(j)}) \, \Delta S_j = g(\bm x^{(i)}), \quad i = 0, 1, \ldots, N,
\]
where $\Delta S_j$ is the area of the grid cell associated with node $\bm y^{(j)}$.

This yields a system of linear equations
\begin{equation} \label{e:6}
A \varphi = f,
\end{equation}
where the right-hand side vector is formed from the observed data:
\[
f_i = g(\bm x^{(i)}), \quad i = 0, 1, \ldots, N.
\]
The matrix $A$ is defined by
\[
A_{ij} = G \frac{x_3^{(i)} - y_3^{(j)}}{|\bm x^{(i)} - \bm y^{(j)}|^3} \, \Delta S_j, \quad i = 0, 1, \ldots, N, \quad j = 0, 1, \ldots, M,
\]
and the unknown vector $\varphi$ corresponds to the discretized surface density:
\[
\varphi_j = \mu(\bm y^{(j)}), \quad j = 0, 1, \ldots, M.
\]

In practice, observational data are subject to measurement noise, typically modeled as normally distributed random errors. Let the perturbed data $\tilde f$ satisfy
\[
\| \tilde f - f \| \le \delta,
\]
where $\delta$ is the noise level. Instead of solving \eqref{e:6} directly, we seek an approximate solution of
\begin{equation} \label{e:7}
A \varphi \approx \tilde f
\end{equation}
in the class of non-negative functions satisfying the constraint
\begin{equation} \label{e:5_repeat}
\varphi \ge 0.
\end{equation}
To solve~\eqref{e:7} under constraint~\eqref{e:5_repeat}, we apply the non-negative least squares method:
\[
\| A \varphi - \tilde f \|^2 \ \longrightarrow \ \min_{\varphi \ge 0}.
\]

The system matrix $A$ in~\eqref{e:6} is bounded but, in general, not symmetric or positive definite. This poses challenges for standard solvers. Common approaches to handling such matrices include iterative methods with symmetrization or regularization techniques.

In our case, the problem is nonlinear due to the a priori constraint on non-negativity. The numerical solution is implemented using the Non-Negative Least Squares algorithm~\cite{lawson1995solving}, which is an adaptation of the classical least squares method using the active set strategy.
To enhance the depth resolution of reconstructed sources and stabilize the solution, regularization techniques based on iterative procedures and noise-aware modeling are recommended.

\section{Numerical experiments}

We consider a model problem with non-dimensionalized parameters, setting the gravitational constant $G = 1$ for simplicity.
Two point sources with masses $0.1$ and $0.2$ are located at coordinates $(-0.2, 0.2, -0.3)$ and $(0.3, -0.1, -0.4)$, respectively.

The observation surface lies at depth $h = 0$ ($x_3 = 0$) and covers a square region $[-1, 1] \times [-1, 1]$.
A regular grid of size $N_1 \times N_2$ is constructed, resulting in a total of $N = (N_1 + 1)(N_2 + 1)$ observation points.
In general, the distribution of points on $\Gamma$ can be arbitrary to account for surface topography, so $x_3 = z(x_1, x_2)$.

\begin{figure}[htb]
\centering
\includegraphics[width=0.45\textwidth]{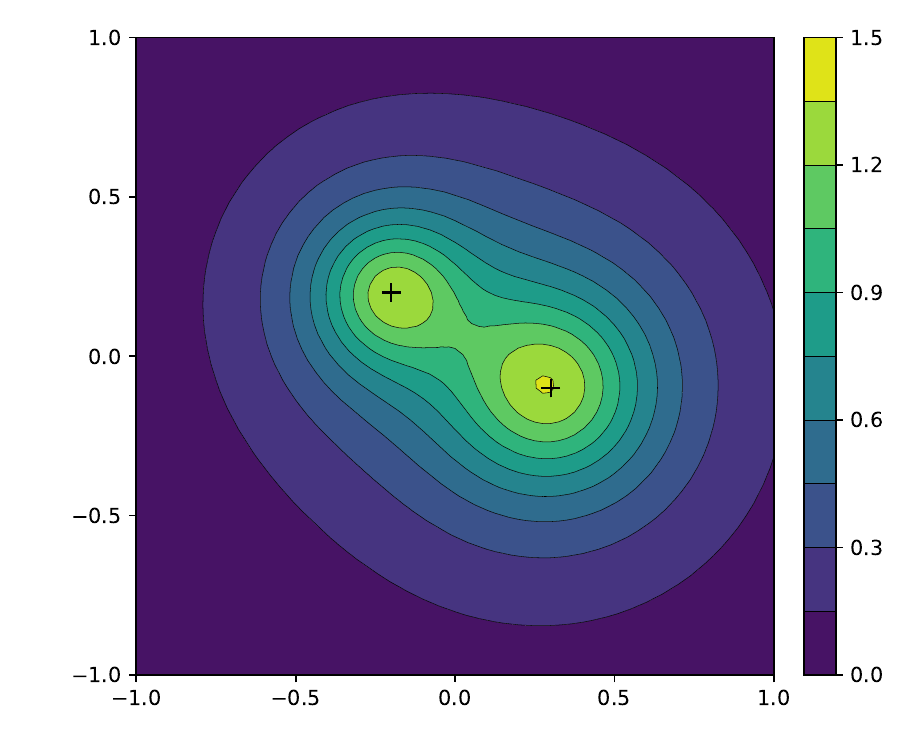}
\includegraphics[width=0.45\textwidth]{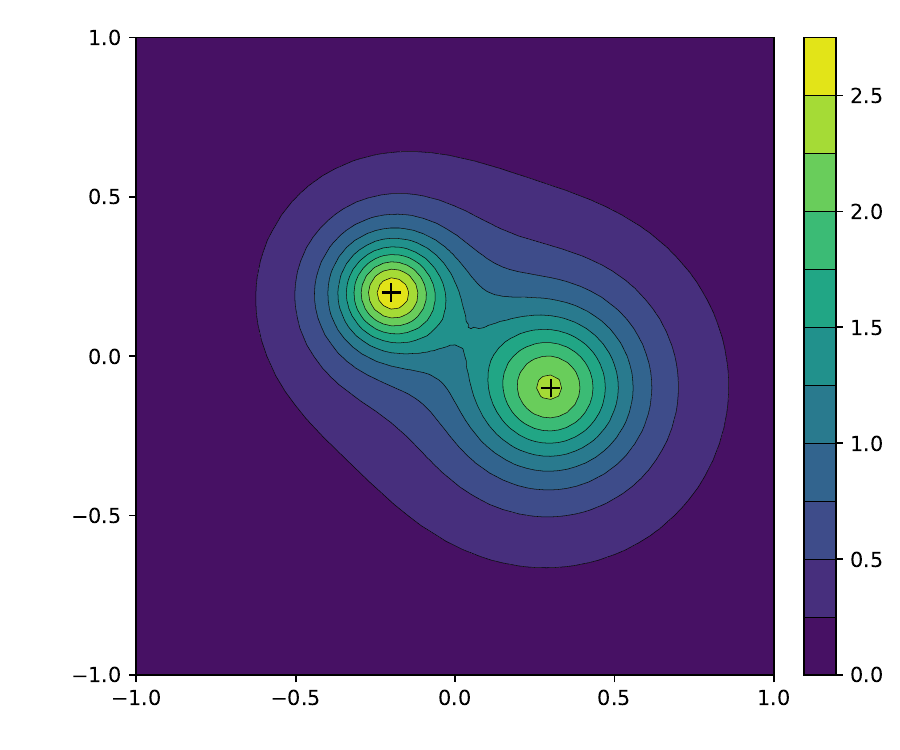}

\includegraphics[width=0.45\textwidth]{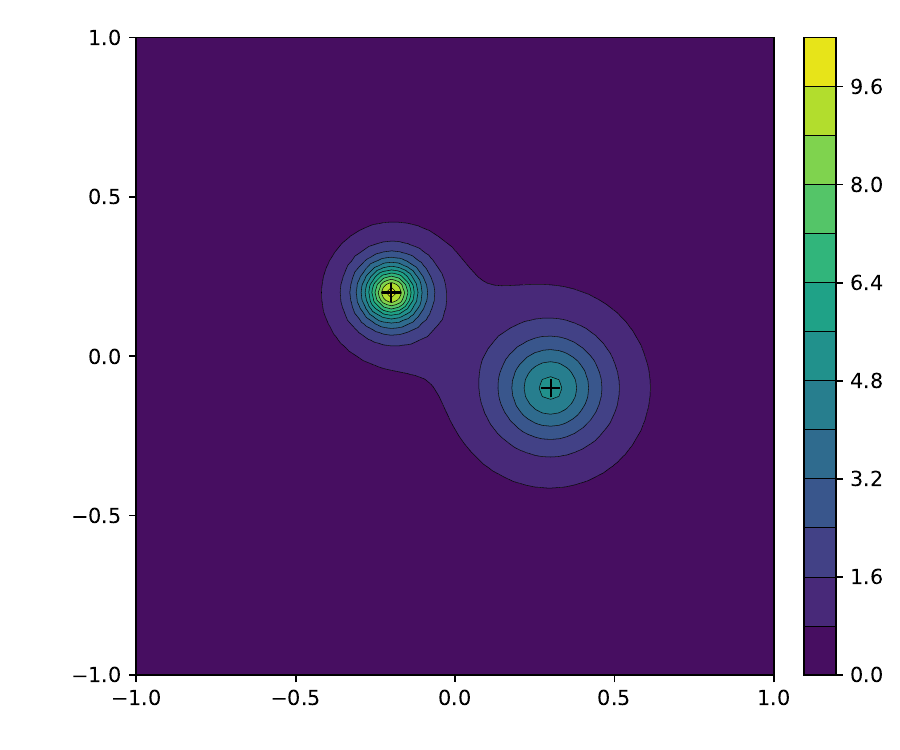}
\includegraphics[width=0.45\textwidth]{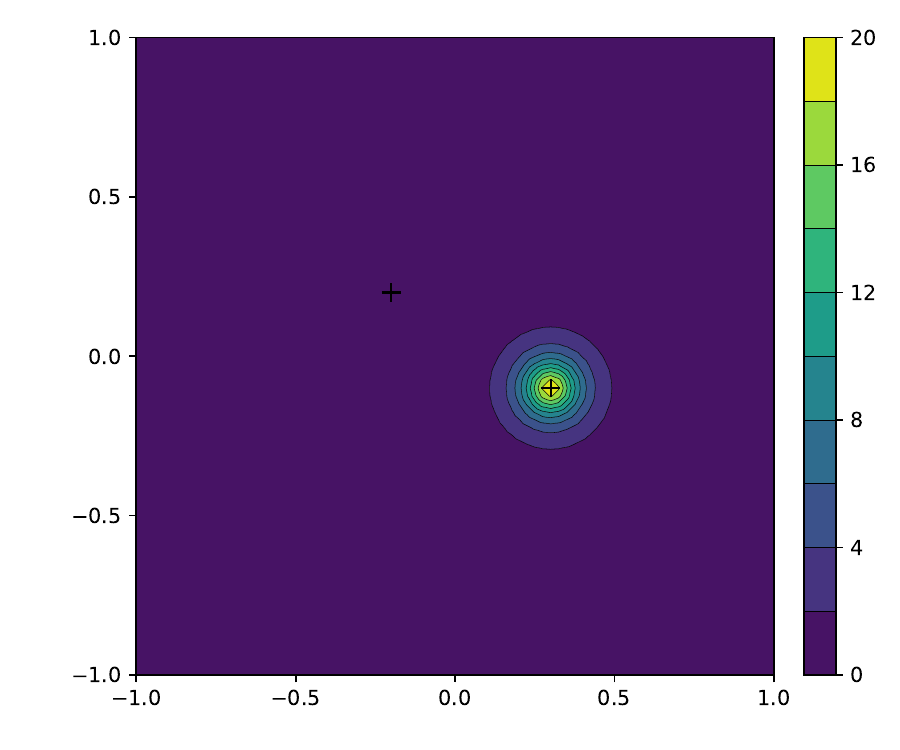}
\caption{Gravitational field $g(\bm x)$ at different depths $h$: from left to right $h=0, 0.1$ (top row), $h=0.2, 0.3$ (bottom row).}
\label{f:501}
\end{figure}

The exact gravitational field at various depths is shown in Figure~\ref{f:501}.
As the depth increases, the localization of the anomalies becomes more pronounced.

We now consider the horizontal continuation surface $\gamma$, a square region $[-1, 1] \times [-1, 1]$, located at depth $h$ ($x_3 = -h$) beneath the observation surface.
A regular grid of size $M_1 \times M_2$ is generated by dividing each side into $M_1$ and $M_2$ intervals, respectively.
The total number of grid nodes (or cells) is $M = (M_1 + 1)(M_2 + 1)$.

\begin{figure}[h!]
\centering
\includegraphics[width=0.45\textwidth]{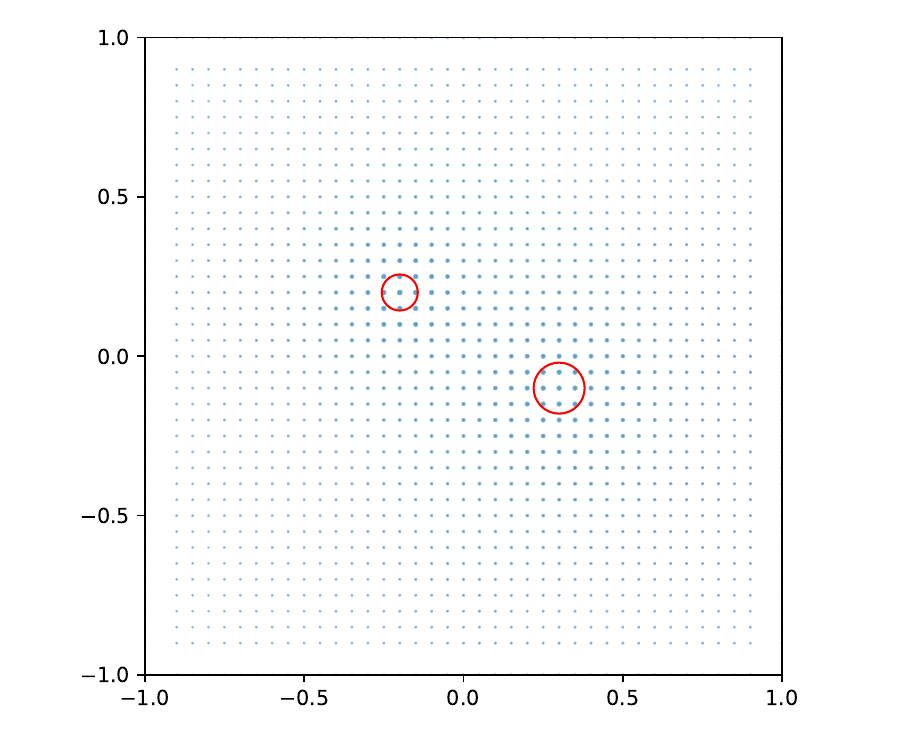}
\includegraphics[width=0.45\textwidth]{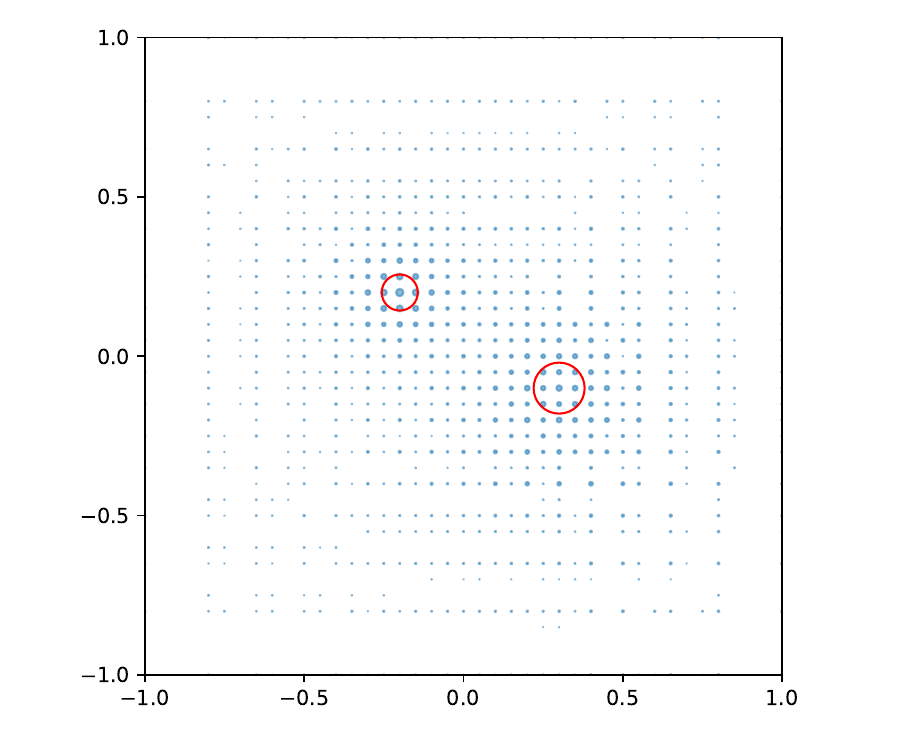}

\includegraphics[width=0.45\textwidth]{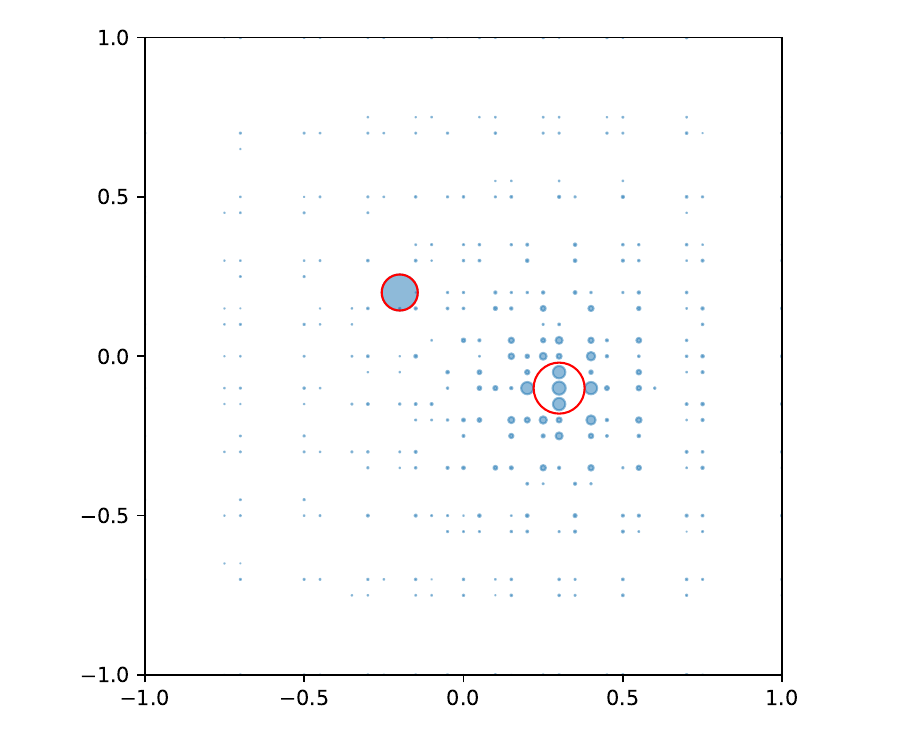}
\includegraphics[width=0.45\textwidth]{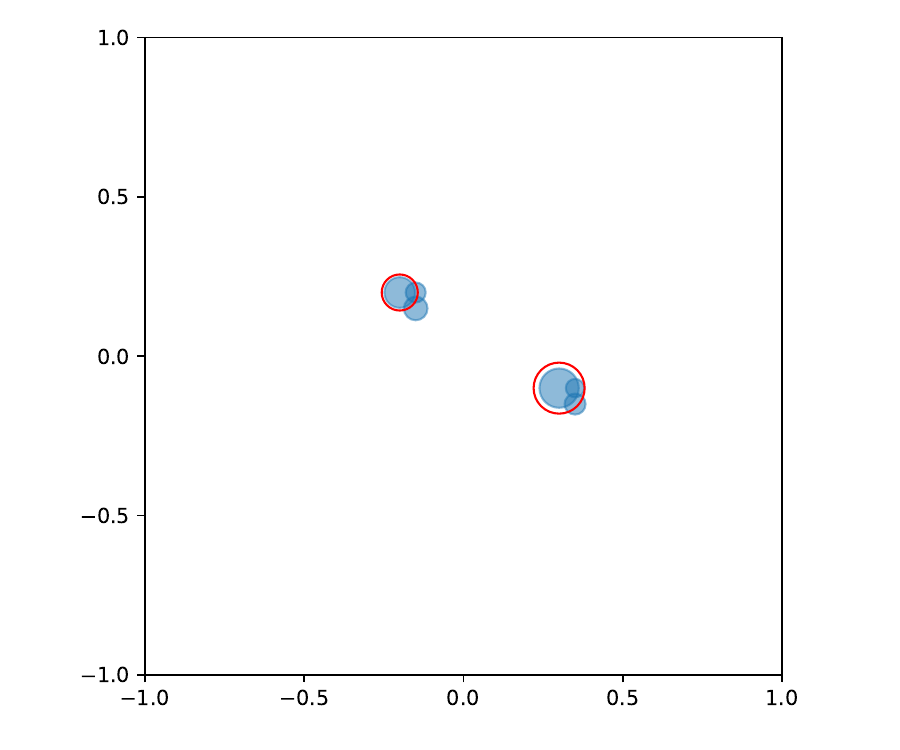}
\caption{Approximate solutions $\mu(\bm x)$ at depths: from left to right $h=0.1, 0.2$ (top row), $h=0.3, 0.4$ (bottom row). The exact point sources are indicated by red circles.}
\label{f:502}
\end{figure}

To perform downward continuation of the gravitational field toward the source masses, a sequence of computations is carried out on a regular $40 \times 40$ grid ($M_1 = M_2 = 40$) over $\gamma$.
Figure~\ref{f:502} shows the reconstructed mass values $\varphi_j \Delta S_j$ at the grid nodes $\bm y^{(j)}$ for different continuation depths.
The true point sources are marked by red circles.
At shallow depths, the reconstructed surface density is spread over the entire continuation domain $\gamma$, and the anomalies are not well localized.
As the depth increases, the solution becomes concentrated near the actual sources, with localization at a few isolated points.

\begin{figure}[htb]
\centering
\includegraphics[width=0.9\textwidth]{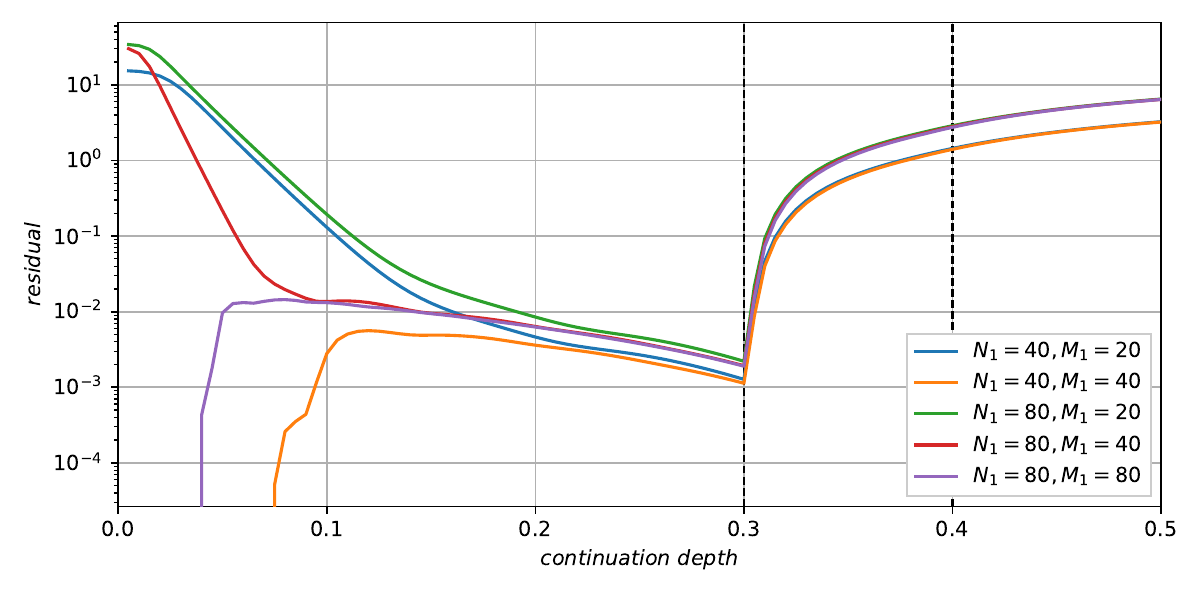}
\caption{Residual $\chi(h)$ at different continuation depths $h$ using noise-free data.}
\label{f:503}
\end{figure}

To quantitatively assess the accuracy of the solution, we compute the residual
\[
\chi(h) = \| A_h \varphi(h) - f \|,
\]
where $\varphi(h)$ is the approximate solution at continuation depth $h$.
The dependence of the residual on depth is shown in Figure~\ref{f:503}, where the depths of the true sources are marked with dashed vertical lines.
When the observation and continuation grids coincide in the $x_1x_2$-plane, the continuation is highly accurate (residual approaches machine precision) at shallow depths.
Otherwise, the residual remains nearly independent of the grid parameters $N_1, N_2$ and $M_1, M_2$.
A sharp decrease in residual, typically corresponding to the depth of the nearest source, indicates a local minimum.

These results suggest a possible algorithm for estimating equivalent source parameters.
First, the depth of the nearest source is estimated from the residual minimum.
Then, the horizontal coordinates and mass are inferred from the reconstructed surface density at that depth.
Subsequent sources are identified by subtracting the contribution of previously identified sources from the input data and repeating the process.

To simulate measurement errors, we introduce noise with relative level $\delta$:
\[
\tilde f_i = f_i + \delta \max_{0 \le k \le N} |f_k| \, \sigma_i, \quad i = 0, 1, \ldots, N,
\]
where $\sigma_i$ are random variables from the standard normal distribution.

\begin{figure}[htb]
\centering
\includegraphics[width=0.9\textwidth]{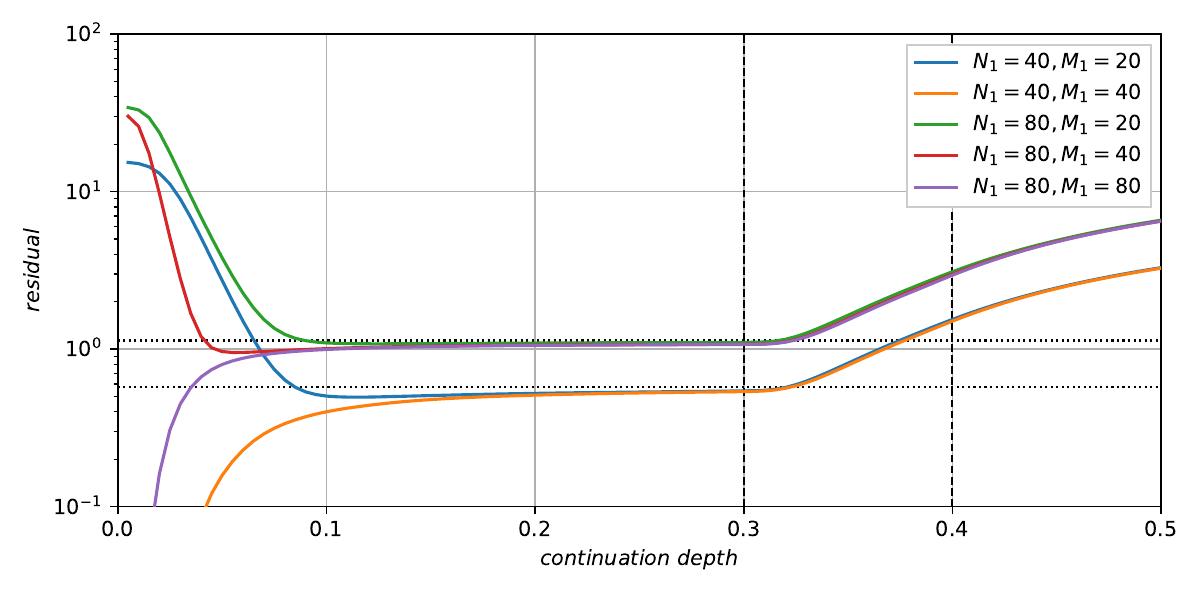}
\caption{Residual $\chi(h)$ at different continuation depths  $h$ with noise level $\delta = 0.01$.}
\label{f:504}
\end{figure}

\begin{figure}[htb]
\centering
\includegraphics[width=0.45\textwidth]{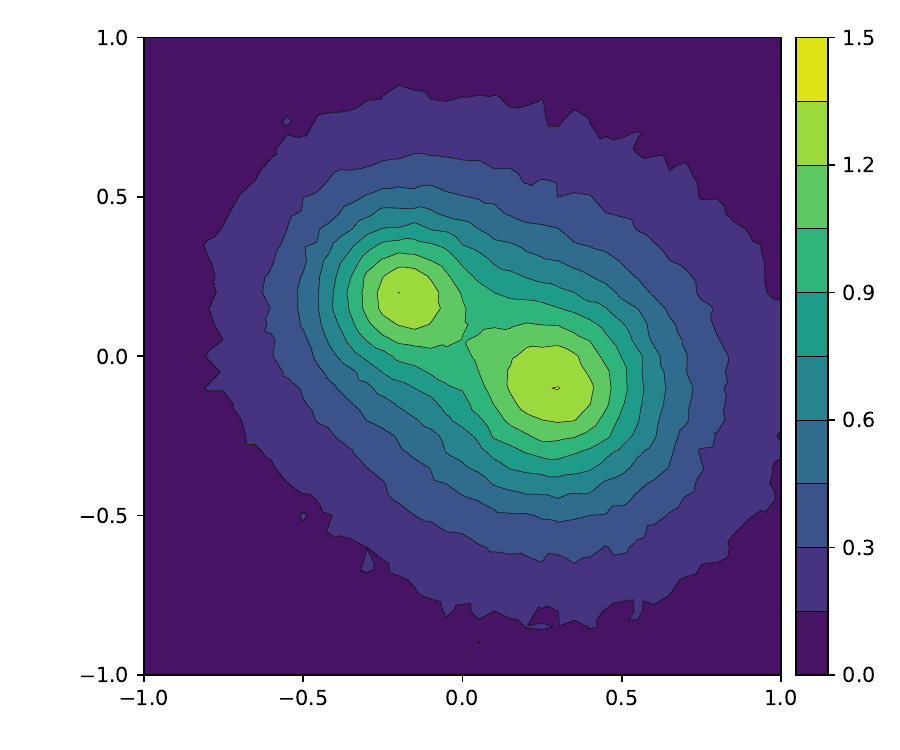}
\includegraphics[width=0.45\textwidth]{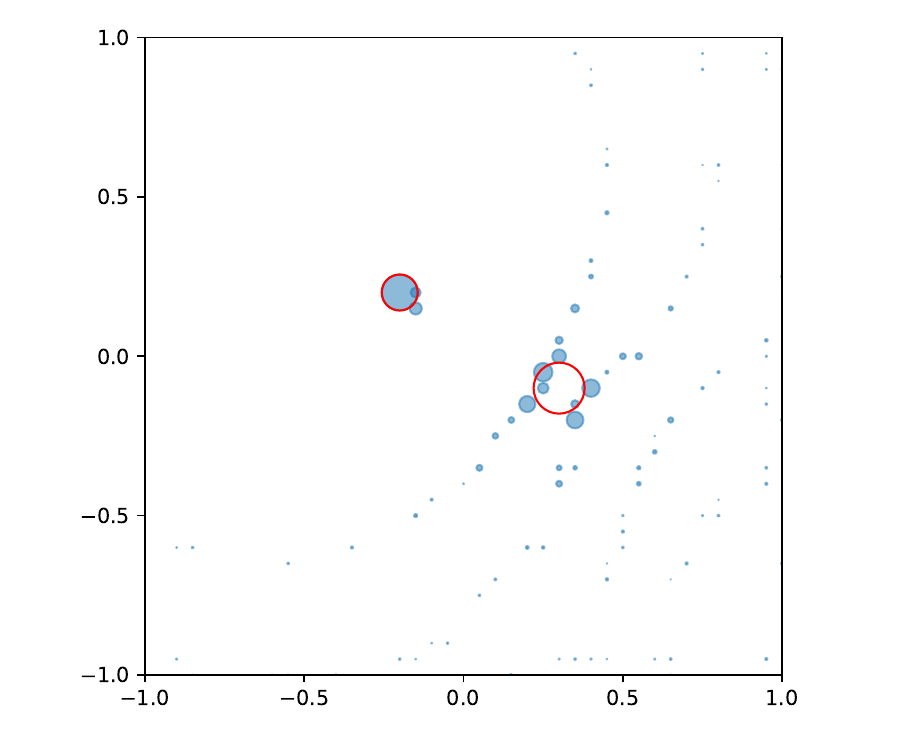}
\caption{Noisy observation data (left) and approximate solution at optimal depth $h^*_{0.01}$ with noise level $\delta = 0.01$ (right).}
\label{f:505}
\end{figure}

The optimal continuation depth is selected as the maximal depth satisfying the residual criterion:
\[
h^*_\delta = \arg\max_h \left\{ \| A_h \varphi(h) - \tilde f \| \le \delta \sqrt{N} \max_i |\tilde f_i| \right\}.
\]
Figure~\ref{f:504} shows the residual profile for $\delta = 0.01$, with the threshold $\delta \sqrt{N} \max_i |\tilde f_i|$ indicated by a dotted horizontal line for different values of $N_1$.
The resulting optimal depth is nearly independent of $N_1$ and $M_1$.
For example, for $N_1 = M_1 = 40$, the optimal depth is $h^*_{0.01} = 0.32$.
The corresponding approximate mass distribution is shown in Figure~\ref{f:505}.

\begin{figure}[htb]
\centering
\includegraphics[width=0.9\textwidth]{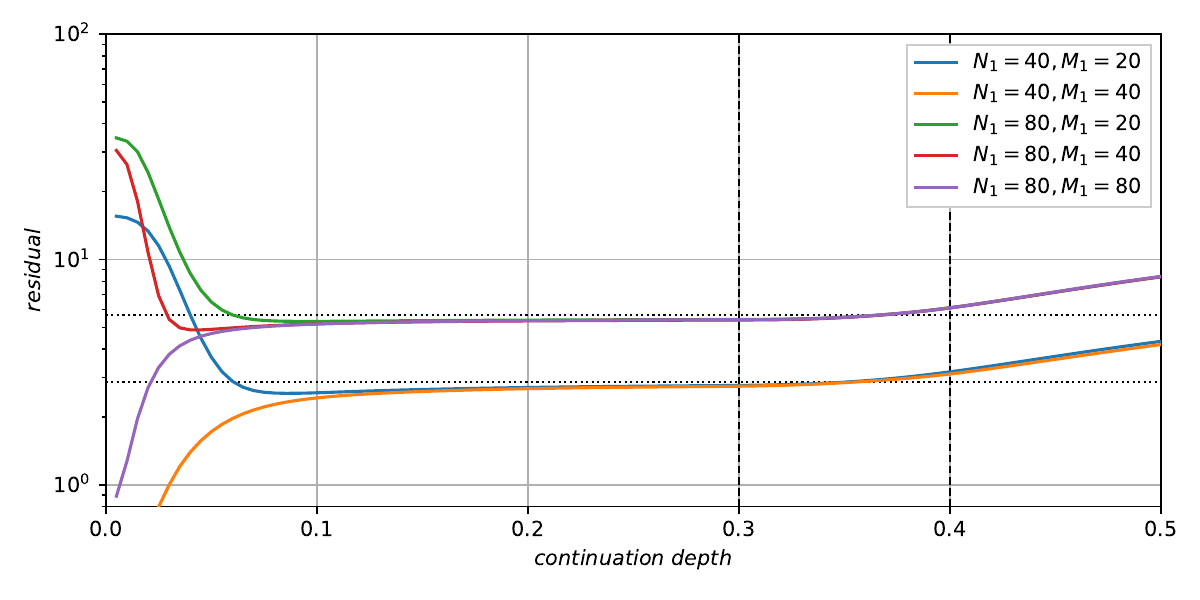}
\caption{Residual $\chi(h)$ at different continuation depths  $h$ with noise level $\delta = 0.05$.}
\label{f:506}
\end{figure}

\begin{figure}[htb]
\centering
\includegraphics[width=0.45\textwidth]{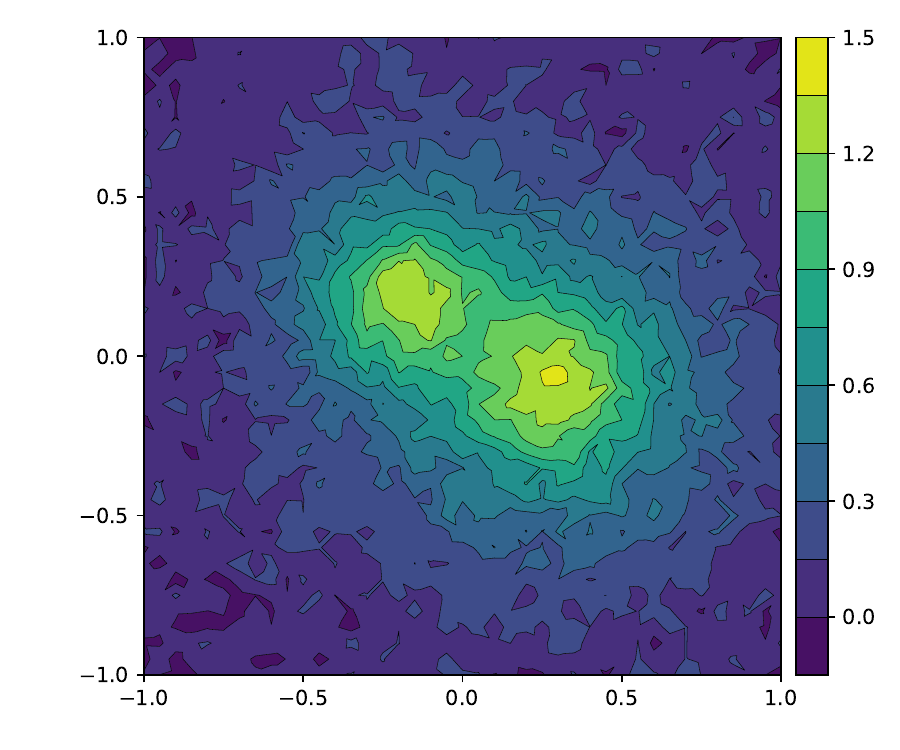}
\includegraphics[width=0.45\textwidth]{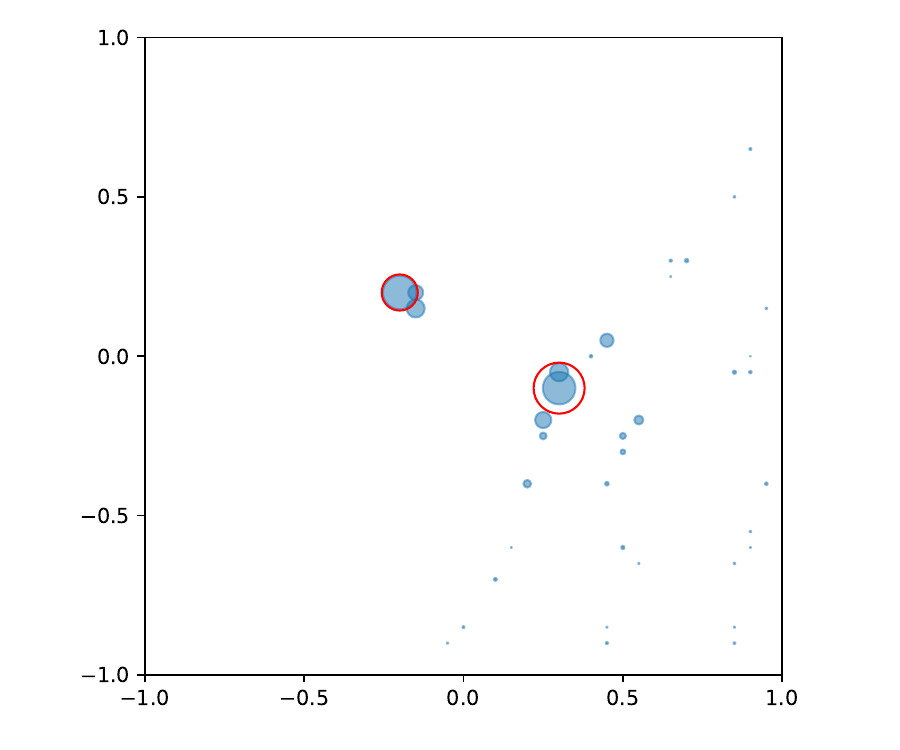}
\caption{Noisy observation data (left) and approximate solution at optimal depth $h^*_{0.05}$ with noise level $\delta = 0.05$ (right).}
\label{f:507}
\end{figure}

Results for $\delta = 0.05$ are shown in Figures~\ref{f:506} and~\ref{f:507}, with the optimal depth $h^*_{0.05} = 0.345$.
The reconstructed solutions remain robust, consistently recovering the source locations within the correct depth range.

\section{Conclusion}

An approximate solution to the ill-posed problem of potential field continuation from the Earth's surface toward underlying sources is obtained using the simple layer potential approach.
This method is widely used in geophysical exploration.

The main results of the study are summarized as follows:
\begin{itemize}
	\item For problems involving sources with constant-sign anomalous density, a class of a priori constraints is identified. The numerical algorithm is based on the non-negative least squares method.
	\item The continuation depth of the simple layer potential plays the role of a regularization parameter, and its selection is guided by the residual principle accounting for measurement noise.
	\item Numerical modeling of three-dimensional continuation problems confirms the stability of the proposed algorithm and its ability to localize sources both horizontally and in depth.
\end{itemize}

\section*{Acknowledgement}

This work was supported by the Ministry of Science and Higher Education of the Russian Federation (Supplementary Agreement No.~075-02-2025-1792, March 11, 2025) and by the Science Committee of the Ministry of Science and Higher Education of the Republic of Kazakhstan (Grant No.~BR27100483).



\begin{thebibliography}{42}
\ifx \bisbn   \undefined \def \bisbn  #1{ISBN #1}\fi
\ifx \binits  \undefined \def \binits#1{#1}\fi
\ifx \bauthor  \undefined \def \bauthor#1{#1}\fi
\ifx \batitle  \undefined \def \batitle#1{#1}\fi
\ifx \bjtitle  \undefined \def \bjtitle#1{#1}\fi
\ifx \bvolume  \undefined \def \bvolume#1{\textbf{#1}}\fi
\ifx \byear  \undefined \def \byear#1{#1}\fi
\ifx \bissue  \undefined \def \bissue#1{#1}\fi
\ifx \bfpage  \undefined \def \bfpage#1{#1}\fi
\ifx \blpage  \undefined \def \blpage #1{#1}\fi
\ifx \burl  \undefined \def \burl#1{\textsf{#1}}\fi
\ifx \doiurl  \undefined \def \doiurl#1{\url{https://doi.org/#1}}\fi
\ifx \betal  \undefined \def \betal{\textit{et al.}}\fi
\ifx \binstitute  \undefined \def \binstitute#1{#1}\fi
\ifx \binstitutionaled  \undefined \def \binstitutionaled#1{#1}\fi
\ifx \bctitle  \undefined \def \bctitle#1{#1}\fi
\ifx \beditor  \undefined \def \beditor#1{#1}\fi
\ifx \bpublisher  \undefined \def \bpublisher#1{#1}\fi
\ifx \bbtitle  \undefined \def \bbtitle#1{#1}\fi
\ifx \bedition  \undefined \def \bedition#1{#1}\fi
\ifx \bseriesno  \undefined \def \bseriesno#1{#1}\fi
\ifx \blocation  \undefined \def \blocation#1{#1}\fi
\ifx \bsertitle  \undefined \def \bsertitle#1{#1}\fi
\ifx \bsnm \undefined \def \bsnm#1{#1}\fi
\ifx \bsuffix \undefined \def \bsuffix#1{#1}\fi
\ifx \bparticle \undefined \def \bparticle#1{#1}\fi
\ifx \barticle \undefined \def \barticle#1{#1}\fi
\bibcommenthead
\ifx \bconfdate \undefined \def \bconfdate #1{#1}\fi
\ifx \botherref \undefined \def \botherref #1{#1}\fi
\ifx \url \undefined \def \url#1{\textsf{#1}}\fi
\ifx \bchapter \undefined \def \bchapter#1{#1}\fi
\ifx \bbook \undefined \def \bbook#1{#1}\fi
\ifx \bcomment \undefined \def \bcomment#1{#1}\fi
\ifx \oauthor \undefined \def \oauthor#1{#1}\fi
\ifx \citeauthoryear \undefined \def \citeauthoryear#1{#1}\fi
\ifx \endbibitem  \undefined \def \endbibitem {}\fi
\ifx \bconflocation  \undefined \def \bconflocation#1{#1}\fi
\ifx \arxivurl  \undefined \def \arxivurl#1{\textsf{#1}}\fi
\csname PreBibitemsHook\endcsname

%%% 1
\bibitem[\protect\citeauthoryear{Gupta}{2011}]{gupta2011encyclopedia}
\begin{bbook}
\bauthor{\bsnm{Gupta}, \binits{H.}}:
\bbtitle{Encyclopedia of Solid Earth Geophysics}.
\bpublisher{Springer},
\blocation{Cham}
(\byear{2011})
\end{bbook}
\endbibitem

%%% 2
\bibitem[\protect\citeauthoryear{Mudrecova and
  Veselov}{1990}]{mudrecova1990gravi}
\begin{bbook}
\beditor{\bsnm{Mudrecova}, \binits{E.A.}},
\beditor{\bsnm{Veselov}, \binits{K.E.}} (eds.):
\bbtitle{Gravimetric Exploration. Geophysicist's Handbook}.
\bpublisher{Nedra},
\blocation{Moscow}
(\byear{1990}).
\bcomment{in {R}ussian}
\end{bbook}
\endbibitem

%%% 3
\bibitem[\protect\citeauthoryear{Parker}{1994}]{Parker1994}
\begin{bbook}
\bauthor{\bsnm{Parker}, \binits{R.L.}}:
\bbtitle{Geophysical Inverse Theory}.
\bpublisher{Princeton University Press},
\blocation{Princeton, NJ}
(\byear{1994})
\end{bbook}
\endbibitem

%%% 4
\bibitem[\protect\citeauthoryear{Blokh}{2009}]{blokh2009interp}
\begin{bbook}
\bauthor{\bsnm{Blokh}, \binits{Y.I.}}:
\bbtitle{Interpretation of Gravity and Magnetic Anomalies}.
\bpublisher{MGGA},
\blocation{Moscow}
(\byear{2009}).
\bcomment{in {R}ussian}
\end{bbook}
\endbibitem

%%% 5
\bibitem[\protect\citeauthoryear{Crossley et~al.}{2013}]{Crossley2013}
\begin{barticle}
\bauthor{\bsnm{Crossley}, \binits{D.}},
\bauthor{\bsnm{Hinderer}, \binits{J.}},
\bauthor{\bsnm{Riccardi}, \binits{U.}}:
\batitle{The measurement of surface gravity}.
\bjtitle{Reports on Progress in Physics}
\bvolume{76}(\bissue{4}),
\bfpage{046101}
(\byear{2013})
\end{barticle}
\endbibitem

%%% 6
\bibitem[\protect\citeauthoryear{Vogel et~al.}{1992}]{Vogel1992}
\begin{bbook}
\beditor{\bsnm{Vogel}, \binits{A.}},
\beditor{\bsnm{Sarwar}, \binits{A.K.M.}},
\beditor{\bsnm{Gorenflo}, \binits{R.}},
\beditor{\bsnm{Kounchev}, \binits{O.I.}} (eds.):
\bbtitle{Theory and Practice of Geophysical Data Inversion}.
\bpublisher{Vieweg+Teubner},
\blocation{Wiesbaden}
(\byear{1992})
\end{bbook}
\endbibitem

%%% 7
\bibitem[\protect\citeauthoryear{Zhdanov}{2002}]{zhdanov2002geophysical}
\begin{bbook}
\bauthor{\bsnm{Zhdanov}, \binits{M.S.}}:
\bbtitle{Geophysical Inverse Theory and Regularization Problems}.
\bpublisher{Elsevier},
\blocation{Amsterdam}
(\byear{2002})
\end{bbook}
\endbibitem

%%% 8
\bibitem[\protect\citeauthoryear{Freeden and
  Nashed}{2018}]{freeden2018handbook}
\begin{bbook}
\beditor{\bsnm{Freeden}, \binits{W.}},
\beditor{\bsnm{Nashed}, \binits{M.Z.}} (eds.):
\bbtitle{Handbook of Mathematical Geodesy: Functional Analytic and Potential
  Theoretic Methods}.
\bpublisher{Birkh{\"a}user}, \blocation{???}
(\byear{2018})
\end{bbook}
\endbibitem

%%% 9
\bibitem[\protect\citeauthoryear{Yagola et~al.}{2014}]{yagola2014obratnye}
\begin{bbook}
\bauthor{\bsnm{Yagola}, \binits{A.G.}},
\bauthor{\bsnm{Yanfei}, \binits{V.}},
\bauthor{\bsnm{Stepanova}, \binits{I.E.}},
\bauthor{\bsnm{Titarenko}, \binits{V.N.}}:
\bbtitle{Inverse Problems and Methods of Their Solution. Applications to
  Geophysics}.
\bpublisher{Binom},
\blocation{Moscow}
(\byear{2014}).
\bcomment{in {R}ussian}
\end{bbook}
\endbibitem

%%% 10
\bibitem[\protect\citeauthoryear{Eppelbaum}{2019}]{eppelbaum2019geophysical}
\begin{bbook}
\bauthor{\bsnm{Eppelbaum}, \binits{L.}}:
\bbtitle{Geophysical Potential Fields: Geological and Environmental
  Applications}.
\bpublisher{Elsevier},
\blocation{Amsterdam}
(\byear{2019})
\end{bbook}
\endbibitem

%%% 11
\bibitem[\protect\citeauthoryear{Blakely}{1996}]{blakely1996potential}
\begin{bbook}
\bauthor{\bsnm{Blakely}, \binits{R.J.}}:
\bbtitle{Potential Theory in Gravity and Magnetic Applications}.
\bpublisher{Cambridge university press},
\blocation{Cambridge}
(\byear{1996})
\end{bbook}
\endbibitem

%%% 12
\bibitem[\protect\citeauthoryear{Dolgal}{2022}]{dolgal2022gravimetry}
\begin{bbook}
\bauthor{\bsnm{Dolgal}, \binits{A.S.}}:
\bbtitle{Gravimetry and Magnetometry: Transformations of Geopotential Fields}.
\bpublisher{Perm State National Research University},
\blocation{Perm}
(\byear{2022}).
\bcomment{in {R}ussian}
\end{bbook}
\endbibitem

%%% 13
\bibitem[\protect\citeauthoryear{Pilkington and
  Boulanger}{2017}]{pilkington2017potential}
\begin{barticle}
\bauthor{\bsnm{Pilkington}, \binits{M.}},
\bauthor{\bsnm{Boulanger}, \binits{O.}}:
\batitle{Potential field continuation between arbitrary surfaces --- comparing
  methods}.
\bjtitle{Geophysics}
\bvolume{82}(\bissue{3}),
\bfpage{9}--\blpage{25}
(\byear{2017})
\end{barticle}
\endbibitem

%%% 14
\bibitem[\protect\citeauthoryear{Pa{\v{s}}teka
  et~al.}{2012}]{pavsteka2012regcont}
\begin{barticle}
\bauthor{\bsnm{Pa{\v{s}}teka}, \binits{R.}},
\bauthor{\bsnm{Karcol}, \binits{R.}},
\bauthor{\bsnm{Ku{\v{s}}nir{\'a}k}, \binits{D.}},
\bauthor{\bsnm{Mojze{\v{s}}}, \binits{A.}}:
\batitle{{REGCONT: A M}atlab based program for stable downward continuation of
  geophysical potential fields using {T}ikhonov regularization}.
\bjtitle{Computers \& Geosciences}
\bvolume{49},
\bfpage{278}--\blpage{289}
(\byear{2012})
\end{barticle}
\endbibitem

%%% 15
\bibitem[\protect\citeauthoryear{Lavrentiev and
  Saveliev}{2006}]{lavrentiev2006operator}
\begin{bbook}
\bauthor{\bsnm{Lavrentiev}, \binits{M.M.}},
\bauthor{\bsnm{Saveliev}, \binits{L.J.}}:
\bbtitle{Operator Theory and Ill-Posed Problems}.
\bpublisher{De Gruyter},
\blocation{Berlin}
(\byear{2006})
\end{bbook}
\endbibitem

%%% 16
\bibitem[\protect\citeauthoryear{Tikhonov and
  Arsenin}{1977}]{tikhonov1977solutions}
\begin{bbook}
\bauthor{\bsnm{Tikhonov}, \binits{A.N.}},
\bauthor{\bsnm{Arsenin}, \binits{V.Y.}}:
\bbtitle{Solutions of Ill-posed Problems}.
\bpublisher{Distributed solely by Halsted Press},
\blocation{New York}
(\byear{1977})
\end{bbook}
\endbibitem

%%% 17
\bibitem[\protect\citeauthoryear{Kabanikhin}{2011}]{kabanikhin2011inverse}
\begin{bbook}
\bauthor{\bsnm{Kabanikhin}, \binits{S.I.}}:
\bbtitle{Inverse and Ill-posed Problems: Theory and Applications}.
\bpublisher{Berlin},
\blocation{De Gruyter}
(\byear{2011})
\end{bbook}
\endbibitem

%%% 18
\bibitem[\protect\citeauthoryear{Samarskii and
  Vabishchevich}{2007}]{samarskii2007numerical}
\begin{bbook}
\bauthor{\bsnm{Samarskii}, \binits{A.A.}},
\bauthor{\bsnm{Vabishchevich}, \binits{P.N.}}:
\bbtitle{Numerical Methods for Solving Inverse Problems of Mathematical
  Physics}.
\bpublisher{De Gruyter},
\blocation{Berlin}
(\byear{2007})
\end{bbook}
\endbibitem

%%% 19
\bibitem[\protect\citeauthoryear{Alifanov et~al.}{1995}]{alifanov1995extrem}
\begin{bbook}
\bauthor{\bsnm{Alifanov}, \binits{O.M.}},
\bauthor{\bsnm{Artyukhin}, \binits{E.A.}},
\bauthor{\bsnm{Rumyantsev}, \binits{S.V.}}:
\bbtitle{Extreme Methods for Solving Ill-Posed Problems with Applications to
  Inverse Heat Transfer Problems}.
\bpublisher{Begell House},
\blocation{New York}
(\byear{1995})
\end{bbook}
\endbibitem

%%% 20
\bibitem[\protect\citeauthoryear{Hansen}{1998}]{hansen1998rank}
\begin{bbook}
\bauthor{\bsnm{Hansen}, \binits{P.C.}}:
\bbtitle{Rank-deficient and Discrete Ill-posed Problems: Numerical Aspects of
  Linear Inversion}.
\bpublisher{SIAM},
\blocation{Philadelphia}
(\byear{1998})
\end{bbook}
\endbibitem

%%% 21
\bibitem[\protect\citeauthoryear{Vasin and Ageev}{2013}]{vasin2013ill}
\begin{bbook}
\bauthor{\bsnm{Vasin}, \binits{V.V.}},
\bauthor{\bsnm{Ageev}, \binits{A.L.}}:
\bbtitle{Ill-Posed Problems with a Priori Information}.
\bpublisher{De Gruyter},
\blocation{Berlin}
(\byear{2013})
\end{bbook}
\endbibitem

%%% 22
\bibitem[\protect\citeauthoryear{Bakushinsky et~al.}{2011}]{Bakushinsky2011}
\begin{bbook}
\bauthor{\bsnm{Bakushinsky}, \binits{A.B.}},
\bauthor{\bsnm{Kokurin}, \binits{M.Y.}},
\bauthor{\bsnm{Smirnova}, \binits{A.}}:
\bbtitle{Iterative Methods for Ill-Posed Problems: An Introduction}.
\bpublisher{De Gruyter},
\blocation{Berlin}
(\byear{2011})
\end{bbook}
\endbibitem

%%% 23
\bibitem[\protect\citeauthoryear{Dolgal et~al.}{2022}]{dolgal2022history}
\begin{barticle}
\bauthor{\bsnm{Dolgal}, \binits{A.S.}},
\bauthor{\bsnm{Pugin}, \binits{A.V.}},
\bauthor{\bsnm{Novikova}, \binits{P.N.}}:
\batitle{History of the method for sourcewise approximations of geopotential
  fields}.
\bjtitle{Izvestiya, Physics of the Solid Earth}
\bvolume{58}(\bissue{2}),
\bfpage{149}--\blpage{171}
(\byear{2022})
\end{barticle}
\endbibitem

%%% 24
\bibitem[\protect\citeauthoryear{Hackbusch}{1995}]{Hackbusch1995}
\begin{bbook}
\bauthor{\bsnm{Hackbusch}, \binits{W.}}:
\bbtitle{Integral Equations: Theory and Numerical Treatment}.
\bpublisher{Birkhäuser},
\blocation{Basel}
(\byear{1995})
\end{bbook}
\endbibitem

%%% 25
\bibitem[\protect\citeauthoryear{Kythe and Puri}{2011}]{kythe2011computational}
\begin{bbook}
\bauthor{\bsnm{Kythe}, \binits{P.}},
\bauthor{\bsnm{Puri}, \binits{P.}}:
\bbtitle{Computational Methods for Linear Integral Equations}.
\bpublisher{Springer},
\blocation{New York}
(\byear{2011})
\end{bbook}
\endbibitem

%%% 26
\bibitem[\protect\citeauthoryear{Tikhonov
  et~al.}{1968}]{tikhonov1968prodolgenie}
\begin{barticle}
\bauthor{\bsnm{Tikhonov}, \binits{A.N.}},
\bauthor{\bsnm{Glasko}},
\bauthor{\bsnm{B.}, \binits{V.}},
\bauthor{\bsnm{Litvinenko}, \binits{O.K.}},
\bauthor{\bsnm{Melikhov}, \binits{V.R.}}:
\batitle{On the downward continuation of the potential toward disturbing masses
  in gravity and magnetic prospecting using the regularization method}.
\bjtitle{Izvestiya, Academy of Sciences of the USSR. Physics of the Solid
  Earth}
\bvolume{2},
\bfpage{30}--\blpage{48}
(\byear{1968})
\end{barticle}
\endbibitem

%%% 27
\bibitem[\protect\citeauthoryear{Tikhonov et~al.}{1995}]{tikhonov1995numerical}
\begin{bbook}
\bauthor{\bsnm{Tikhonov}, \binits{A.N.}},
\bauthor{\bsnm{Goncharsky}, \binits{A.V.}},
\bauthor{\bsnm{Stepanov}, \binits{V.V.}},
\bauthor{\bsnm{Yagola}, \binits{A.G.}}:
\bbtitle{Numerical Methods for the Solution of Ill-Posed Problems}.
\bpublisher{Springer},
\blocation{Dordrecht}
(\byear{1995})
\end{bbook}
\endbibitem

%%% 28
\bibitem[\protect\citeauthoryear{Strakhov and Stepanova}{2002}]{strakhov2002s}
\begin{barticle}
\bauthor{\bsnm{Strakhov}, \binits{V.N.}},
\bauthor{\bsnm{Stepanova}, \binits{I.E.}}:
\batitle{The {S}-approximation method and its application to gravity problems}.
\bjtitle{Izvestiya, Physics of the Solid Earth}
\bvolume{38}(\bissue{2}),
\bfpage{91}--\blpage{107}
(\byear{2002})
\end{barticle}
\endbibitem

%%% 29
\bibitem[\protect\citeauthoryear{Stepanova et~al.}{2019}]{stepanova2019approx}
\begin{barticle}
\bauthor{\bsnm{Stepanova}, \binits{I.E.}},
\bauthor{\bsnm{Kerimov}, \binits{I.A.}},
\bauthor{\bsnm{Yagola}, \binits{A.G.}}:
\batitle{Approximation approach in various modifications of the method of
  linear integral representations}.
\bjtitle{Izvestiya, Physics of the Solid Earth}
\bvolume{55}(\bissue{2}),
\bfpage{218}--\blpage{231}
(\byear{2019})
\end{barticle}
\endbibitem

%%% 30
\bibitem[\protect\citeauthoryear{Stepanova
  et~al.}{2020}]{stepanova2020separation}
\begin{barticle}
\bauthor{\bsnm{Stepanova}, \binits{I.E.}},
\bauthor{\bsnm{Raevsky}, \binits{D.N.}},
\bauthor{\bsnm{Shchepetilov}, \binits{A.V.}}:
\batitle{Separation of potential fields generated by different sources based on
  modified {S}-approximations}.
\bjtitle{Izvestiya, Physics of the Solid Earth}
\bvolume{56}(\bissue{3}),
\bfpage{379}--\blpage{391}
(\byear{2020})
\end{barticle}
\endbibitem

%%% 31
\bibitem[\protect\citeauthoryear{Li and Oldenburg}{2010}]{li2010rapid}
\begin{barticle}
\bauthor{\bsnm{Li}, \binits{Y.}},
\bauthor{\bsnm{Oldenburg}, \binits{D.W.}}:
\batitle{Rapid construction of equivalent sources using wavelets}.
\bjtitle{Geophysics}
\bvolume{75}(\bissue{3}),
\bfpage{51}--\blpage{59}
(\byear{2010})
\end{barticle}
\endbibitem

%%% 32
\bibitem[\protect\citeauthoryear{Balk et~al.}{2016}]{balk2016effective}
\begin{barticle}
\bauthor{\bsnm{Balk}, \binits{P.I.}},
\bauthor{\bsnm{Dolgal}, \binits{A.S.}},
\bauthor{\bsnm{Pugin}, \binits{A.V.}},
\bauthor{\bsnm{Michurin}, \binits{A.V.}},
\bauthor{\bsnm{Simanov}, \binits{A.A.}},
\bauthor{\bsnm{Sharkhimullin}, \binits{A.F.}}:
\batitle{Effective algorithms for sourcewise approximation of geopotential
  fields}.
\bjtitle{Izvestiya, Physics of the Solid Earth}
\bvolume{52}(\bissue{6}),
\bfpage{896}--\blpage{911}
(\byear{2016})
\end{barticle}
\endbibitem

%%% 33
\bibitem[\protect\citeauthoryear{Temirbekov
  et~al.}{2022}]{temirbekov2022numerical}
\begin{barticle}
\bauthor{\bsnm{Temirbekov}, \binits{N.}},
\bauthor{\bsnm{Temirbekova}, \binits{L.}},
\bauthor{\bsnm{Nurmangaliyeva}, \binits{M.}}:
\batitle{Numerical solution of the first kind fredholm integral equations by
  projection methods with wavelets as the basis functions}.
\bjtitle{TWMS Journal of Pure and Applied Mathematics}
\bvolume{13}(\bissue{1}),
\bfpage{105}--\blpage{118}
(\byear{2022})
\end{barticle}
\endbibitem

%%% 34
\bibitem[\protect\citeauthoryear{Dampney}{1969}]{dampney1969equivalent}
\begin{barticle}
\bauthor{\bsnm{Dampney}, \binits{C.N.G.}}:
\batitle{The equivalent source technique}.
\bjtitle{Geophysics}
\bvolume{34}(\bissue{1}),
\bfpage{39}--\blpage{53}
(\byear{1969})
\end{barticle}
\endbibitem

%%% 35
\bibitem[\protect\citeauthoryear{Leonov et~al.}{2024}]{leonov2024solving}
\begin{barticle}
\bauthor{\bsnm{Leonov}, \binits{A.S.}},
\bauthor{\bsnm{Lukyanenko}, \binits{D.V.}},
\bauthor{\bsnm{Yagola}, \binits{A.G.}}:
\batitle{Solving some inverse problems of gravimetry and magnetometry using an
  algorithm that improves matrix conditioning}.
\bjtitle{Computational Mathematics and Mathematical Physics}
\bvolume{64}(\bissue{10}),
\bfpage{2178}--\blpage{2193}
(\byear{2024})
\end{barticle}
\endbibitem

%%% 36
\bibitem[\protect\citeauthoryear{Zeng et~al.}{2013}]{zeng2013adaptive}
\begin{barticle}
\bauthor{\bsnm{Zeng}, \binits{X.}},
\bauthor{\bsnm{Li}, \binits{X.}},
\bauthor{\bsnm{Su}, \binits{J.}},
\bauthor{\bsnm{Liu}, \binits{D.}},
\bauthor{\bsnm{Zou}, \binits{H.}}:
\batitle{An adaptive iterative method for downward continuation of
  potential-field data from a horizontal plane}.
\bjtitle{Geophysics}
\bvolume{78}(\bissue{4}),
\bfpage{43}--\blpage{52}
(\byear{2013})
\end{barticle}
\endbibitem

%%% 37
\bibitem[\protect\citeauthoryear{Li et~al.}{2023}]{li2023stable}
\begin{barticle}
\bauthor{\bsnm{Li}, \binits{D.}},
\bauthor{\bsnm{Liu}, \binits{X.}},
\bauthor{\bsnm{Zhai}, \binits{Z.}},
\bauthor{\bsnm{Du}, \binits{J.}},
\bauthor{\bsnm{Chen}, \binits{C.}},
\bauthor{\bsnm{Ma}, \binits{J.}},
\bauthor{\bsnm{Wang}, \binits{Y.}}:
\batitle{A stable downward continuation method for processing gravity data
  using the equivalent sources with compactness and smoothing constraints}.
\bjtitle{Journal of Applied Geophysics}
\bvolume{215},
\bfpage{105128}
(\byear{2023})
\end{barticle}
\endbibitem

%%% 38
\bibitem[\protect\citeauthoryear{Vabishchevich}{2024}]{vabishchevich2024computational}
\begin{barticle}
\bauthor{\bsnm{Vabishchevich}, \binits{P.N.}}:
\batitle{Computational algorithm for the continuation of potential fields
  towards gravitating masses [in {R}ussian]}.
\bjtitle{Numerical Methods and Programming}
\bvolume{25}(\bissue{1}),
\bfpage{1}--\blpage{9}
(\byear{2024})
\end{barticle}
\endbibitem

%%% 39
\bibitem[\protect\citeauthoryear{Lawson and Hanson}{1995}]{lawson1995solving}
\begin{bbook}
\bauthor{\bsnm{Lawson}, \binits{C.L.}},
\bauthor{\bsnm{Hanson}, \binits{R.J.}}:
\bbtitle{Solving Least Squares Problems}.
\bpublisher{SIAM},
\blocation{New {Y}ork}
(\byear{1995})
\end{bbook}
\endbibitem

%%% 40
\bibitem[\protect\citeauthoryear{Latt{\`e}s and Lions}{1969}]{lattes1969method}
\begin{bbook}
\bauthor{\bsnm{Latt{\`e}s}, \binits{R.}},
\bauthor{\bsnm{Lions}, \binits{J.L.}}:
\bbtitle{The Method of Quasi-Reversibility: Applications to Partial
  Differential Equations}.
\bpublisher{American Elsevier Publishing Company},
\blocation{New York}
(\byear{1969})
\end{bbook}
\endbibitem

%%% 41
\bibitem[\protect\citeauthoryear{Vabishchevich and Pulatov}{2002}]{vab2002num}
\begin{barticle}
\bauthor{\bsnm{Vabishchevich}, \binits{P.N.}},
\bauthor{\bsnm{Pulatov}, \binits{P.A.}}:
\batitle{Numerical solution of a problem of potential field continuation}.
\bjtitle{Matematicheskoe Modelirovanie}
\bvolume{14}(\bissue{6}),
\bfpage{91}--\blpage{104}
(\byear{2002}).
\bcomment{in {R}ussian}
\end{barticle}
\endbibitem

%%% 42
\bibitem[\protect\citeauthoryear{Vabishchevich and
  Pulatov}{1983}]{vabishchevich1983economical}
\begin{botherref}
\oauthor{\bsnm{Vabishchevich}, \binits{P.N.}},
\oauthor{\bsnm{Pulatov}, \binits{P.A.}}:
Economical difference methods for solving direct problems in gravimetric and
  magnetic prospecting.
Izv. Akad. Nauk SSSR: Fiz. Zemli,
(10),
68--76
(1983).
in {R}ussian
\end{botherref}
\endbibitem

\end{thebibliography}
\end{document}